\documentclass[letterpaper, 10 pt, conference]{ieeeconf}
\IEEEoverridecommandlockouts
\usepackage[absolute]{textpos}
\usepackage{amsmath}
\usepackage{amssymb}
\usepackage{graphicx}
\usepackage{subfig}
\usepackage{algorithm}
\usepackage{booktabs}
\title{\LARGE \bf Fast Trajectory Planning for Automated Vehicles using \mbox{Gradient-based Nonlinear Model Predictive Control}}
\author{Franz Gritschneder, Knut Graichen and Klaus Dietmayer%
\thanks{The authors are with Institute of Measurement, Control and Microtechnology, Ulm University, 89081 Ulm, Germany
        {\mbox{(e-mail:\tt\small firstname.lastname@uni-ulm.de)}}}%
}
\begin{document}
\begin{textblock*}{170mm}(17.5mm,9.3375mm)
\parindent0mm
\normalfont\scriptsize
Copyright $\copyright$ 2018 IEEE.
Personal use of this material is permitted.
Permission from IEEE must be obtained for all other uses, in any current or future media, including reprinting/republishing this material for advertising or promotional purposes, creating new collective works, for resale or redistribution to servers or lists, or reuse of any copyrighted component of this work in other works.
\end{textblock*}
\maketitle
\thispagestyle{empty}
\pagestyle{empty}
\begin{abstract}
Motion trajectory planning is one crucial aspect for automated vehicles, as it governs the own future behavior in a dynamically changing environment.
A good utilization of a vehicle's characteristics requires the consideration of the nonlinear system dynamics within the optimization problem to be solved.
In particular, real-time feasibility is essential for automated driving, in order to account for the fast changing surrounding, e.g. for moving objects.
The key contributions of this paper are the presentation of a fast optimization algorithm for trajectory planning including the nonlinear system model.
Further, a new concurrent operation scheme for two optimization algorithms is derived and investigated.
The proposed algorithm operates in the submillisecond range on a standard PC.
As an exemplary scenario, the task of driving along a challenging reference course is demonstrated.
\end{abstract}
\section{INTRODUCTION}
Fully automated driving is still an unsolved problem.
Trajectory planning is one key enabling module for automated vehicles, as its outcome describes the vehicle's future motion.
The physical behavior of the motion in time is subject to constraints originating from the vehicle's environment as well as comfort and safety goals.
The first aspect ensures a comfortable motion profile for the passengers, whereas rapid changes in acceleration have to be avoided.
Second, the motion has to be collision free anyway, regarding other static as well as dynamic objects.

In general, the task of automated driving is even more complicated.
While static objects can be considered easily, the future behavior of dynamic objects is hard to predict as their behavior can change abruptly.
For example, a pedestrian can alter his direction or speed in an unprecedented manner.
Consistently, this illuminates the need for a motion planning module with a fast computation performance as well as scalability for the planning horizon.

Trajectory planning for automated vehicles is an active field of research for more than a decade.
First approaches \cite{Kelly2003},\cite{Howard2010} already used a model predictive control formulation for vehicle applications.
A first approach \cite{Werling2010} exploited the decomposition of the planning task into a lateral and longitudinal motion profile search.
Recent work \cite{Schildbach2016a}, \cite{Yi2016a}, \cite{Gutjahr2017} successfully implemented model predictive control (MPC) algorithms by a combined consideration of a vehicle's nonlinear dynamics.
The nonholonomic character makes the task of trajectory planning even more complicated.
However, this approach characterizes the real system behavior best but raises difficulties in computing a real-time solution for the optimal control problem.
The approaches from above concern the solution of an optimal control problem by first discretize and then optimize.

\section{CONTRIBUTIONS}
In this paper, a first optimize and then discretize solution strategy is presented.
The continuous formulation of the nonlinear optimal control problem covers the real dynamics best.
Further, two instances of the solution approach are run in a hierarchical fashion, where one computes a reference course for the second instance, which calculates the actual motion profile for a vehicle.
The main contribution of our work addresses the above mentioned challenges and is twofold.
The solution of the optimization problem uses a nonlinear model predictive controller.
\mbox{(i) Our} approach employs a fast real-time capable solution method.
(ii) A hierarchical operation scheme using two MPC instances, which are run in parallel, is introduced.
The latter is beneficial if different domains for the prediction horizon are present.

The remainder of the paper is structured as follows.
First, the vehicle model together with the nonlinear MPC formulation and solution strategy is introduced (\ref{sec:vmpc}).
Next (\ref{sec:pmpc}), a thorough discussion for the course model follows.
In the final part (\ref{sec:impl}), the concurrent operation scheme algorithm is established and experimental results are presented.
\section{MPC for vehicle dynamics}
\label{sec:vmpc}
The desired task of automated driving relies on an appropriate choice of model especially for a model predictive control (MPC) approach.
This section introduces the nonlinear vehicle model, and formulates the cost function for the MPC scheme.
Thereafter, the solution strategy of the MPC formulation is derived.
\subsection{Vehicle model}
\begin{figure}[t]
\centering
\includegraphics[page=1]{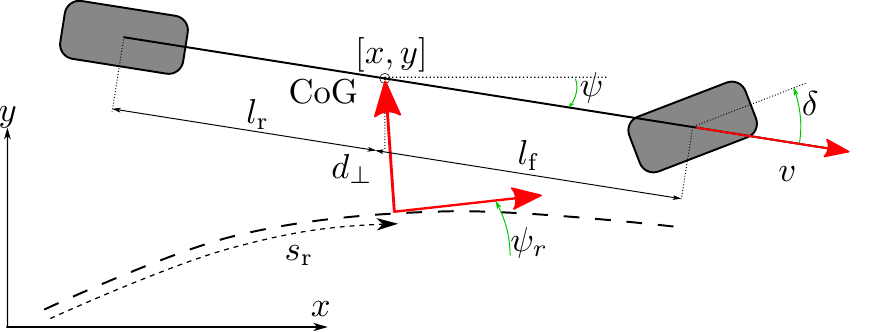}
\caption{Sketch of the vehicle model with its system states together with the reference curve (dashed line).}
\label{fig:vehiclemodel}
\end{figure}

Figure \ref{fig:vehiclemodel} sketches the kinematic model of a vehicle.
The model originates from the bicycle model and was introduced by \cite{Werling2012}, having the state vector defined as
\begin{align}
\boldsymbol{x}_{\text{v}} = \begin{bmatrix} {x}, {y}, {\psi}, {\delta}, {v}, {d_{\perp}}, {\psi_{r}}, {s_{\text{r}}} \end{bmatrix}^{\intercal}.
\end{align}
The system's state vector $\boldsymbol{x}_{\text{v}}$ constitutes the position $({x},{y})$ of the vehicle's center of gravity (CoG) in a cartesian coordinate system, the yaw angle ${\psi}$, the steer angle ${\delta}$ of the front wheels and the velocity ${v}$.
The remaining three state variables $({d_{\perp}}, {\psi_{r}}, {s_{\text{r}}})$ cover the relation of the vehicle's center of gravity with respect to the reference course.
Figure \ref{fig:vehiclemodel} represents the reference course as dashed line, where ${d_{\perp}}$ is the signed lateral displacement of the vehicle from the reference course, ${\psi_{r}}$ is the reference orientation, or in other words the tangent vector of the reference course, and ${s_{\text{r}}}$ traces the driven arc length along the reference course.
The system's input $\boldsymbol{u}_{v}$ reads two entities, namely the steer rate ${\dot{\delta}}$ and the acceleration $\dot{{v}}$.
The first-order differential equations
\begin{align}
{\dot{x}}     ~=~ & {v} \cdot \cos({\psi})\\
{\dot{y}}     ~=~ & {v} \cdot \sin({\psi})\\
{\dot{\psi}} ~=~ & \frac{{v}}{ \ell \cdot \left( 1 + (\frac{{v}}{v_{ch}})^2  \right) } \cdot {\delta}\\
{\dot{\delta}} ~=~ & u_{\text{v},1}\\
{\dot{v}}     ~=~ & u_{\text{v},2}\\
{\dot{d}_{\perp}}~=~ & {v} \cdot \sin({\psi} - {\psi_{r}})\label{eq:refeq1}\\
{\dot{\psi}_{r}}~=~ & {v} \cdot \frac{\cos({\psi} - {\psi_{r}})}{1 - {d_{\perp}} \cdot \kappa({s_{\text{r}}})} \kappa({s_{\text{r}}})\label{eq:refeq2}\\
{\dot{s}_{\text{r}}}~=~ & {v} \cdot \frac{\cos({\psi} - {\psi_{r}})}{1 - {d_{\perp}} \cdot \kappa({s_{\text{r}}})},\label{eq:refeq3}
\end{align}
describe the nonlinear vehicle dynamics $\dot{\boldsymbol{x}}_{\text{v}} = f(\boldsymbol{x}_{\text{v}}, \boldsymbol{u}_{\text{v}})$.
Particularly, the vehicle's dynamics is dependent on two parameters, the length $\ell$ and the so called characteristic velocity $v_{\text{ch}}$ (\cite{Werling2012}).
The first describes the overall length of the vehicle, while the second merges the length from the center of gravity to the front respectively to the rear, the tire stiffness from the front and rear tires together with the vehicle's mass.
The formula is given by $ v_{\text{ch}} = \sqrt{\frac{\ell^2 c_f c_r}{m\left(c_r l_r -c_f l_f\right)}} $.
Instead of computing a geometrical projection onto a reference line,
equations (\ref{eq:refeq1})-(\ref{eq:refeq3}) express the projection of the vehicle's center of gravity onto the reference course by means of differential equations already.
The functional $\kappa({s_{\text{r}}})$ describes the curvature of the reference course along the arc length ${s_{\text{r}}}$.
Up to now, we assume $\kappa({s_{\text{r}}})$ for the vehicle dynamical model is given and fixed.
\subsection{MPC formulation}
A general MPC formulation is presented, together with an efficient online solver approach using gradient based optimization.
Afterwards the application of the previously introduced nonlinear vehicle model is discussed.\\
The employed MPC scheme relies on a continuous optimal control problem (OCP) formulated as
\begin{align}
&\min_{ \boldsymbol u(\cdot) } && J(\boldsymbol  u) = \int_{0}^{T_{\text{hor}}} \! l_{\text{v}}(\boldsymbol x(t),\boldsymbol u(t), t) \,  {\rm d} t\label{eq:mpcformulation}\\
&\text{s.t.} && \dot{\boldsymbol{x}_{\text{v}}}(t) = \boldsymbol f(\boldsymbol x, \boldsymbol u, t), \quad \boldsymbol x(t_i) = \boldsymbol x_0 \label{eq:xvdyn}\\
&&& \boldsymbol h(\boldsymbol x(t),\boldsymbol u(t),t) \le \boldsymbol 0 \quad \forall t\in[t_i,t_f]\,\label{eq:mpcsicons}
\end{align}
with $T_{\text{hor}}$ as the prediction horizon, $\dot{\boldsymbol{x}_{\text{v}}}$ describing the system dynamics (\ref{eq:xvdyn}) having the initial state $\boldsymbol{x}_0$.
State constraints are handled by the inequality equation (\ref{eq:mpcsicons}).
\subsection{Cost function and constraints for the vehicle model}
Applying the vehicle model to the MPC optimization scheme, the cost function and state as well as input constraints need to be defined.
First, the state $\boldsymbol{\xi}$ as \mbox{$\mathbb{R}^8 \to \mathbb{R}^6$}, \mbox{$\boldsymbol x(t) \mapsto \boldsymbol \xi(t)$} is defined.
This mapping describes a transformation from system states to representative states as
\begin{align}
\boldsymbol{\xi} =& \Big[ \begin{matrix}v_{\perp}&a_{\perp}&{\psi}&\kappa_{v}&{d_{\perp}}&v\end{matrix} \Big]^{\intercal}\\
         =& \Big[ \begin{matrix}v {\psi} &v^2 {\psi}&{\psi}&\kappa_{v}&{d_{\perp}}&v\end{matrix}\Big]^{\intercal}
\end{align}
with $\kappa_{v} = \delta\cdot\left(\ell\left(1+\left(\frac{v}{v_{\text{ch}}}\right)^{2}\right)\right)^{-1}$.
Using this state representation, the desired behavior of vehicle motion is then expressed within the cost function.
The first two states of $\boldsymbol{\xi}$ describe the lateral velocity respectively the lateral acceleration.
The third and fourth state give the orientation and curvature value.
Other than \cite{Werling2012}, where a collision avoidance scenario was presented and the lateral absolute deviation from the reference line as well as the velocity have not been taken into account in the cost terms ultimately, our task of trajectory planning must incorporate these two entities necessarily.
Thus, the function $l_v(\boldsymbol{x}_{\text{v}},\boldsymbol{u}_{\text{v}})$ terms the integral costs according to (\ref{eq:mpcformulation}) and is defined as
\begin{align}
\begin{aligned}
l_v( \boldsymbol x, \boldsymbol u) = &~ (\boldsymbol{\xi}-\boldsymbol{\hat{\xi}})^T \boldsymbol{Q} (\boldsymbol{\xi}-\boldsymbol{\hat{\xi}})
  +   \boldsymbol{u}_{\text{v}}^T \boldsymbol{R} \boldsymbol{u}_{\text{v}}\\
    + &~ \gamma*p({\delta}, {\delta}^{-}, {\delta}^{+}),
\end{aligned}
\end{align}
penalizing deviations from their respective desired values $\boldsymbol{\hat{\xi}}$. The cost for the input states $\boldsymbol{u}_{\text{v}}$ are computed based on the absolute state values.
In order to account for the state constraints
\begin{align}
h_1(\boldsymbol x) = &~ \delta - \delta_{\text{max}}  \leq 0\label{eq:h1vmpc}\\
h_2(\boldsymbol x) = &~ \delta + \delta_{\text{min}}  \geq 0\label{eq:h2vmpc}
\end{align}
a soft-penalty function is introduced, ensuring valid steer angles ${\delta}$.
The penalty function $p(x,x^{-},x^{+})$ is given by:
\begin{align}
p(x,x^{-}, x^{+}) \mathrel{\vcenter{\baselineskip0.5ex \lineskiplimit0pt \hbox{\scriptsize.}\hbox{\scriptsize.}}}=
\begin{cases}
      (x-x^{+})^2 & \text{if } x~\text{\textgreater}~x^{+},\\
      (x-x^{-})^2 & \text{if } x~\text{\textless}~x^{-},\\
      0 & \text{otherwise}.
    \end{cases}\label{eq:penfunc}
\end{align}
Finally, the inequality constraint
\begin{align}
h_3(\boldsymbol x) = &~ v^2\cdot \left|\left|\kappa\right|\right| \leq 0,\label{eq:h3vmpc}
\end{align}
expresses an absolute maximum lateral acceleration.
Clearly, the desired behavior features an anticipated reduction of velocity right ahead of a curve, which is encoded within the target value of the velocity $\hat{v}_{\text{tar}}$.
Given the curvature value $\kappa$ and an upper absolute lateral acceleration value $a_{\perp,\text{max}}$, the maximum velocity calculates as 
\begin{align}
v_{\text{max}} =  \sqrt{\frac{\left|\left|a_{\perp,\text{max}}\right|\right|}{\left|\left|\kappa_{v}\right|\right|} }.\label{eq:maxalatformula}
\end{align}
Accounting for speed limit values that are typically below the maximum velocity due to the lateral acceleration constraint on straight roads, the minimum function combines the speed limit value $v_{\text{sl}}$ with the value from (\ref{eq:maxalatformula}) using the minimum operator $\hat{v}_{\text{tar}} =  \min( v_{\text{sl}}, v_{\text{max}} ).$
Finally, the detailed definition of the target state $\boldsymbol{\hat{\xi}}$ is given by
\begin{align}
\boldsymbol{\hat{\xi}}=& \begin{bmatrix} \hat{v}_{\perp} & \hat{a}_{\perp} & {\psi} & \hat{\kappa} & \hat{{d_{\perp}}} & \hat{v} \end{bmatrix}^{\intercal}\\
=& \begin{bmatrix} \hat{v}_{\text{tar}}{\psi_{r}}&\,
  \hat{v}_{\text{tar}}^2\kappa_r&
{\psi_{r}}&
\kappa_r&
\hat{{d_{\perp}}}&
\hat{v}_{\text{tar}}
\end{bmatrix}^{\intercal},
\end{align}
with $r$ indices denoting the reference course values.
The remaining constraints for the control states, that is
\begin{align}
h_4(\boldsymbol{u}_{\text{v}}) = &\,{\dot{\delta}} - {\dot{\delta}}_{\text{min}} &\!\geq 0 ~\wedge~ h_5(\boldsymbol{u}_{\text{v}}) =&\,{\dot{\delta}} + {\dot{\delta}}_{\text{max}}\!\!& \leq 0\label{eq:h45vmpc}\\
h_6(\boldsymbol{u}_{\text{v}}) = &\,a - a_{\text{min}}         &\!\geq 0 ~\wedge~ h_7(\boldsymbol{u}_{\text{v}}) =&\,a + a_{\text{max}}\! & \leq 0\label{eq:h67vmpc}
\end{align}
are considered in the optimization algorithm, presented now.
\subsection{Gradient based optimization strategy}
Accounting for the real-time requirements in automated driving, we use a computationally efficient solver, based on a gradient method.
The gradient method exploits the first order optimality conditions.
The solution strategy uses the Hamiltonian function
\begin{align}
\begin{aligned}
H(\tau, \boldsymbol x(\tau), \boldsymbol(\tau), \boldsymbol \lambda(\tau)) = & l(\boldsymbol x(\tau), \boldsymbol u(\tau), \tau)~+\\
&\boldsymbol \lambda^{T}(\tau) f(\boldsymbol x(\tau), \boldsymbol u(\tau), \tau +t_k)
\end{aligned}
\end{align}
where $\boldsymbol{\lambda}$ describes the vector of costates.
The first order optimality conditions then follow from Pontryagins Maximum Principle \cite{Kirk1970}, \cite{Berkovitz1974},
\begin{align}
\dot{x}^{*}_{k}(\tau)       = & f(\tau+t_k, x^{*}_{k}(\tau), u^{*}_{k}(\tau)), ~~~~~~~~~~ x^{*}_{k} = x_k\label{eq:optcon1}\\
\dot{\lambda}^{*}_{k}(\tau) = & -H_{x}(\tau,x^{*}_{k}(\tau), u^{*}_{k}(\tau)), ~~~~~~~~~~ \lambda^{*}_{k} = 0\label{eq:optcon2}\\
\dot{u}^{*}_{k}(\tau)       = & \underset{u\in [u^{-}, u^{+}]}{\arg \min} H(\tau, x^{*}_{k}(\tau), u, \lambda^{*}_{k}(\tau)), ~~ \tau \in [0, T]\label{eq:optcon3}
\end{align}
where $H_{\boldsymbol{x}}$ shortcuts the partial derivative $\partial H/\partial \boldsymbol{x}$.
The optimality conditions (\ref{eq:optcon1}) - (\ref{eq:optcon3}) are iteratively solved using a gradient method for the minimization problem in (\ref{eq:mpcformulation}), outlined in algorithm \ref{alg:gradalg}.
\begin{algorithm}[ht]
\caption{Gradient Algorithm}
\label{alg:gradalg}
\begin{itemize}
\renewcommand{\labelitemi}{\scriptsize$\blacksquare$}
\item Initialization of input trajectory \mbox{$\boldsymbol u^{(0)} (\tau) \in \left[u^{-}, u^{+}\right]$} $\forall~\tau \in \left[0,T\right]$.
\item Gradient iterations for $j=0,...,M$:
\begin{enumerate}
\item[1) ] Forward integration of the system dynamics with initial state $\boldsymbol x_{k}^{(j)} = x_k$, that is
\begin{align}
\boldsymbol{\dot{x}}_{k}^{(j)} = \boldsymbol f(\boldsymbol x_{k}^{(j)}, \boldsymbol u_{k}^{(j)})\label{eq:fwdInt}
\end{align}
\item[2) ] Backward integration of the adjoint dynamics with the terminal state $\boldsymbol \lambda_{k}^{(j)} = \boldsymbol{0}$
\begin{align}
\boldsymbol{\dot{\lambda}}_{k}^{(j)} = \nabla_{\boldsymbol x} H(\boldsymbol x_{k}^{(j)}, \boldsymbol u_{k}^{(j)}, \lambda_{k}^{(j)})
\end{align}
\item[3) ] Evaluation of the gradient $\boldsymbol{g}_{k}^{(j)}$ with
\begin{align}
\boldsymbol g_{k}^{(j)} = \nabla_{\boldsymbol u} H(\boldsymbol x_{k}^{(j)}, \boldsymbol u_{k}^{(j)}, \boldsymbol \lambda_{k}^{(j)})
\end{align}
\item[4) ] (Approximate) step size computation as a line search problem by
\begin{align}
\alpha_{k}^{(j)} = \arg \underset{\alpha > 0}{\min}~J\left( \Psi(\boldsymbol u_{k}^{(j)} - \alpha \boldsymbol g_{k}^{(j)}), \boldsymbol x_{k} \right),\label{eq:lsp}
\end{align}
with search direction $\boldsymbol{g}_k^{(j)}$ projection function
\begin{align}
\Psi(\boldsymbol u) \mathrel{\vcenter{\baselineskip0.5ex \lineskiplimit0pt \hbox{\scriptsize.}\hbox{\scriptsize.}}}=
\begin{cases}
\boldsymbol u^{-} & \text{if } \boldsymbol u\, \text{\textless}\, \boldsymbol u^{-}\\
\boldsymbol u^{+} & \text{if } \boldsymbol u\, \text{\textgreater}\, \boldsymbol u^{+}\\
\boldsymbol u & \text{otherwise}
\end{cases}\label{eq:incon}
\end{align}
\item[5)] Control update \mbox{$\boldsymbol u_{k}^{(j)}(\tau) = \Psi \left( \boldsymbol u_{k}^{(j)}-\alpha^{i} \boldsymbol g_{k}^{(j)} \right),$} \mbox{$~~\tau \in \left[0, T \right]$}.
\end{enumerate}
\end{itemize}
\end{algorithm}
The actual minimization of the cost function happens by a gradient descent step.
Since the step size is not computable in an analytic fashion, the step size is calculated through an approximate solution for equation (\ref{eq:lsp}) using the formula from \cite{Barzilai1988}
\begin{align}
\alpha_{k}^{(j)} = \frac{\int_{0}^{T_h} \left(\boldsymbol g_{k}^{(j)} - \boldsymbol g_{k}^{(j-1)}\right)^{\intercal} \left(\boldsymbol u_{k}^{(j)} - \boldsymbol u_{k}^{(j-1)}\right) {\rm d}\tau}
                        {\int_{0}^{T_h} \left(\boldsymbol g_{k}^{(j)} - \boldsymbol g_{k}^{(j)-1}\right)^{\intercal} \left(\boldsymbol g_{k}^{(j)} - \boldsymbol g_{k}^{(j)-1}\right) {\rm d}\tau}.
\end{align}
This formula reduces the computational burden for the line search problem in (\ref{eq:lsp}) significantly preserving good results at the same time.
Further, the input constraints are handled by equation (\ref{eq:incon}), based on a projection function.
\section{MPC for curvature approximation}
\label{sec:pmpc}
So far, it was assumed that the curvature value $\kappa$ was given.
In the next lines, the retrieval of the functional $\kappa({s_{\text{r}}})$, used in equation (\ref{eq:refeq2}) and (\ref{eq:refeq3}) is discussed.
Analogous to the MPC formulation for the vehicle dynamics, we formulate a second MPC problem for the task of reference course approximation.
\subsection{Reference course model}
The uniquely describing property for curves is the curvature value $\kappa$ \cite{Pressley2001}.
Let $\mathbb{R} \to \mathbb{R}^2,~ t \in [0,1], t \mapsto \boldsymbol{p}(t)$ with $\boldsymbol{p}(t) = \begin{bmatrix}x(t)&y(t)\end{bmatrix}^{\intercal}$ describe a parameterized 2D curve, with points $p(t)$ and argument $t$.
The formula for computing the curvature\cite{Pressley2001}, where the dot notation refers to the derivation along the function argument $t$, is
\begin{align}
\kappa(t) = & \frac{\dot{x}\ddot{y} - \ddot{x}\dot{y}}{\sqrt{(\dot{x}^2 + \dot{y}^2)}^{3}}\label{kappaFormula}.
\end{align}
However, the required functional $\kappa({s_{\text{r}}})$ according to equation (\ref{eq:refeq2}) and (\ref{eq:refeq3}), is a parameterized functional along the arc length value $s$.
Now, a possible strategy for the transformation of the parameterizing argument is outlined.
The curvature $\kappa(s)$ at point $\boldsymbol{p}(t)$ is defined as the rate of change of the tangent vector $\varphi$.
Using the chain rule, one can write $\kappa(s) = \frac{\partial \varphi(t)}{\partial t}\frac{\partial s}{\partial t}$.
Thus, the course length function $s(t)$, parameterized by $t$ is necessary for the transformation procedure.
In general, a closed form solution for the line integral $s(t)$ for the functional $\boldsymbol{p}(t)$ does not exist.
Practically, in order to circumvent this problem, numerical formulas for the equation (\ref{kappaFormula}) facilitate the approximate computation for $\kappa$ as $\tilde{\kappa}$, using sampled points from $\boldsymbol{p}(t)$.
Performing a numerical integration using the computed curvature values $\tilde{\kappa}$ results in a significant error.

\begin{figure}[t]
\centering
\includegraphics[page=2]{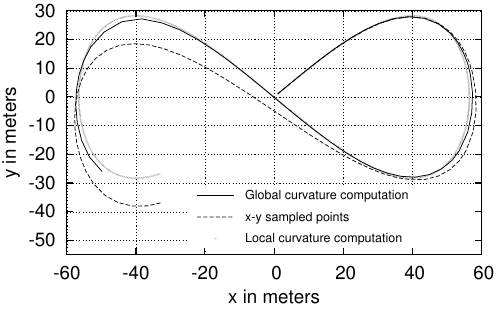}
\caption{Lying eight reference curve, given by the 2D-sample points in x- and y- coordinates. The dashed line arises after integration of the computed curvature values $\kappa$ using formula (\ref{kappaFormula}) for each sampled point. Applying the global optimization approach via the MPC formulation from (\ref{eq:cmpcformulation}) returns the black solid line.}
\label{fig:lyingEightCurve}
\end{figure}

Figure \ref{fig:lyingEightCurve} illustrates the error propagation result by the dashed line, after numerical integration of the curvature values $\tilde{\kappa}$.
Clearly, the curve is observable distinct from the desired course.
This undesirable phenomenon becomes especially critical in real-world traffic scenarios, where errors in lateral direction of the reference line have dramatic consequences.
Starting from sampled points as input data representation, that is $p_i \in \mathbb{R}^2$, \mbox{$i \in [0, N-1]$} with $N \in \mathbb{N}$ denoting the number of sampled points,
an optimization based approach for the computation of the curvature values $\kappa(s)$ is proposed.
In consequence, the aforementioned errors are alleviated considerably.
\subsection{MPC formulation}
Similar to (\ref{eq:mpcformulation}), the MPC problem formulates as
\begin{align}
&\min_{ \boldsymbol u(\cdot) } && J(\boldsymbol  u) =\int_{0}^{S_{\text{hor}}} \! l_{\text{c}}(\boldsymbol{x}_{\text{c}}(s),\boldsymbol u_{\text{c}}(s)) \,  {\rm d} s\label{eq:cmpcformulation}\\
&\text{s.t.} && \boldsymbol{{x}'}_{\text{c}}(s) = \boldsymbol f_{\text{c}}(\boldsymbol{x}_{\text{c}}, \boldsymbol u_{\text{c}}), \quad \boldsymbol{x}_{\text{c}}(s_i) = \boldsymbol{x}_{\text{c},0},
\end{align}
with $\boldsymbol{x}_{\text{c}}$ denoting the course state vector and $\boldsymbol{u}_{\text{c}}$ the control input.
A distinct property compared to (\ref{eq:mpcformulation}) is, that the prediction horizon domain now becomes the course length value $s$.
Accordingly, the dash-notation ${(\cdot)}' = \frac{\partial (\cdot)}{\partial s}$ for the system dynamic
\begin{align}
\boldsymbol{{x}'}_{\text{c}} = \boldsymbol{f}_{\text{c}}(\boldsymbol{x}_{\text{c}}, \boldsymbol u_{\text{c}}) = 
\begin{bmatrix} x^{\prime}\\y^{\prime}\\\varphi^{\prime}\end{bmatrix} =
\begin{bmatrix}\cos(\varphi)\\\sin(\varphi)\\\kappa\end{bmatrix}
\end{align}
is used.
Moreover, constraints for state and inputs are not necessary.
The system dynamic for the reference course model describes, how the curvature $\kappa$ governs the reference course.
The $[x,y]^{\intercal} $ coordinates describe the 2D position of a curve point at the arc length $s$, while the angle $\varphi$ is the tangent vector.
The control input for this dynamical system description is a scalar value, the curvature value $\kappa$.
Next, the cost function for the reference course model is outlined in the following lines.
\subsection{Cost function}
The cost function $l_{\text{c}}(\boldsymbol{x}_{\text{c}}(s),\boldsymbol u_{\text{c}}(s))$ concerns the error for the curve point positions with respect to the given sampled points by a squared euclidean distance measure.
Also, the input control is penalized, prohibiting overfitting effects.
However, since an MPC scheme based on a continuous formulation is applied, the sampled discrete points need to be integrated to the continuous cost function appropriately.
Therefore, we compute a linear interpolation in the fashion of divided differences \cite{DivDiffPaper} for the evaluation of $\hat{x}(s)$ and $\hat{y}(s)$ using the discretized points.
The sampled points $p_i$ are augmented by their course length value $s_i$, that is $\tilde{p}_i = [x_i,y_i,s_i]^{\intercal}$.
For $s\in\left[0,S_{\text{h}}\right]$ and $s_j$ and the N sampled points $p_i = [x_i,y_i]^{\intercal},~i\in[0,N-1]$,  the numerically robust formula is then given by
\begin{align}
\hat{x}(s) = x_{j-1}\frac{s_j-s}{s_j-s_{j-1}} + x_j\frac{s-s_{j-1}}{s_j-s_{j-1}},
\end{align}
where $s_{j-1} \leq s \leq s_j$ holds, and analogous computes $\hat{y}(s)$.
Finally, the cost function is defined as
\begin{align}
\begin{aligned}
l_{\text{c}}(\boldsymbol{x}_{\text{c}}, \boldsymbol u_{\text{c}}(s)) =& \begin{bmatrix}x(s) - \hat{x}(s)\\y(s) - \hat{y}(s)\end{bmatrix}^{\intercal}
\boldsymbol{Q}_{\text{c}}
\begin{bmatrix}x(s) - \hat{x}(s)\\y(s) - \hat{y}(s)\end{bmatrix} +\\
&\boldsymbol u_{\text{c}}^{\intercal}\boldsymbol{R_{\text{c}}}\boldsymbol u_{\text{c}},
\end{aligned}
\end{align}
with $\boldsymbol{Q}_{\text{c}}$ and $\boldsymbol{R}_{\text{c}}$ being positive definite diagonal matrices.
\section{Implementation and Results}
\label{sec:impl}
\begin{figure*}[t]
  \centering
  \subfloat[Velocity trajectory with target value as dashed gray underlay.
            The overshoot at second 20 occurs on the straight path where the lateral displacement has more importance for the optimization.
            The control trajectory $u_{\text{v},2} $, that is the longitudinal acceleration hits the upper bound at 18 seconds once.]
            {\resizebox{86mm}{!}{\includegraphics[page=3]{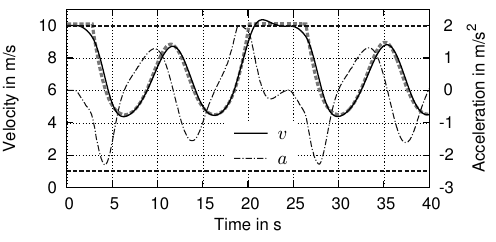}}\label{fig:f1}}
  \hfill
  \subfloat[Steer angle trajectory and steer rate as control input.]
            {\resizebox{86mm}{!}{\includegraphics[page=5]{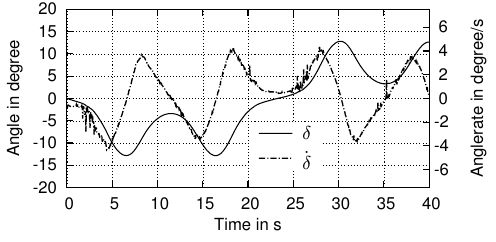}}\label{fig:f2}}
  \vskip\baselineskip
  \subfloat[Lateral offset to the reference track.
            This entity is especially important for automated driving and clearly within desired bounds $(< 1\text{cm})$.]
            {\resizebox{86mm}{!}{\includegraphics[page=4]{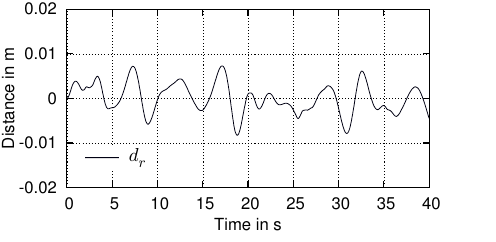}}\label{fig:f3}}
  \hfill
  \subfloat[Orientation trajectory of the reference curve and the vehicle.]
           {\resizebox{86mm}{!}{\includegraphics[page=6]{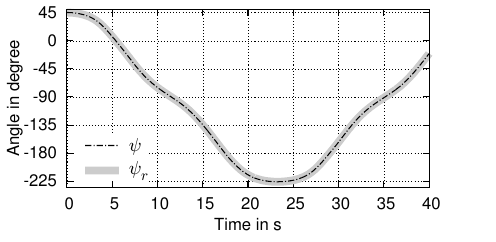}}\label{fig:f4}}
  \caption{State and control trajectory graphs over time for the entire drive along the lying eight of Figure \ref{fig:lyingEightCurve}.}
  \label{fig:stategraphs}
\end{figure*}
This section first describes a concurrent operation scheme using two MPC instances.
One MPC instance holds a vehicle dynamics model, as described in \mbox{Section (\ref{sec:vmpc})} while the second MPC addresses the course approximation problem.
Demonstrating the performance of the proposed MPC solution scheme, experimental results are discussed in the second part of this section.
\subsection{Concurrent MPC instances}
\begin{figure}[ht]
\centering
\includegraphics[page=9]{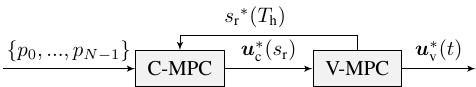}
\caption{The block diagram illustrates the interconnection of the two concurrent running MPC instances.}
\label{fig:blockdiagram}
\end{figure}

The concurrent solution strategy as illustrated in \mbox{Figure \ref{fig:blockdiagram}}, is devised.
Given the sampled 2D-points $p_i$ and an initial state $x_0$ for the vehicle system, the solution scheme uses \emph{two} MPC instances simultaneously in a concurrent operation mode.
The first MPC instance, which we name C-MPC (\emph{C} for \emph{C}ourse), uses the 2D-points and computes a curvature trajectory $u^{*} = \kappa(s)$ for an initial horizon $S_0$ (see \mbox{Fig. \ref{fig:blockdiagram}}).
Carrying out a numerical integration for the optimal control trajectory of the \mbox{C-MPC}, that is $\boldsymbol{u}^{*}_{\text{c}}$, Figure \ref{fig:lyingEightCurve} draws the resulting curve exemplarily.
The C-MPC forwards its solution, that is $\boldsymbol{u}^{*}_{\text{c}}$ to the \emph{V-MPC} (\emph{V} for \emph{V}ehicle), that computes the actual motion trajectory for the vehicle.
The optimization horizon of the \mbox{C-MPC} is dependent on the prediction horizon of the \mbox{V-MPC}.
Thus, the \mbox{V-MPC} hands the prediction horizon variable ${s_{\text{r}}}^{*}(T)$ back to the C-MPC instance.
The block diagram in Figure \ref{fig:blockdiagram} visualizes the data flow, while algorithm \ref{alg:twompcalg} gives the algorithmic steps.
[ht]
The simulation environment was set up within MATLAB.
A vanilla version of the GRAMPC toolbox \cite{Kaepernick2014ECC} provides the solver routines in plain C code.
\begin{algorithm}[th]
\caption{Concurrent MPC operation algorithm}
\label{alg:twompcalg}
\begin{itemize}
\renewcommand{\labelitemi}{\scriptsize$\blacksquare$}
 \item Initialize C-MPC and V-MPC with $S_{\text{h}}=S_0$ and $\boldsymbol{\boldsymbol{x}_{\text{v}}}_{\text{0}}$
 \item Run both MPC instances in parallel as follows:
 \begin{enumerate}
   \item[1) ] Execute gradient algorithm \ref{alg:gradalg} for each MPC instance, that is C-MPC and V-MPC in \emph{parallel}.
   \item[2) ] Mutual update:
   \begin{itemize}
     \item[a) ] Update the C-MPC integration horizon $S$ using the V-MPC state ${s_{\text{r}}}^{*}(T_{\text{h}})$ from equation (\ref{eq:fwdInt}).
     \item[b) ] Update the $\kappa({s_{\text{r}}})$ with the control trajectory from the C-MPC, that is $\boldsymbol{u}_{\text{c}}^{*}(s_{\text{r}})$.
   \end{itemize}
 \end{enumerate}
\end{itemize}
\end{algorithm}
We conducted experiments for driving along straight lines, constant curvature tracks, that is circles, as well as left and right turns.
The most challenging scenario observed was the lying eight which we present now.
Figure \ref{fig:lyingEightCurve} visualizes the chosen reference course.
The asymmetric curvature course ranges from zero to \mbox{$\pm$ 0.075 $1/\text{m}$}, that is an approximated circle with 13 meters of radius.
Together with a straight part at the center of the curve, this altogether covers a representative variety of state combinations.
\subsection{Parameter design}
\begin{figure*}[t]
  \centering
  \subfloat[Acceleration after a sharp right turn.]
           {\resizebox{86mm}{!}{\includegraphics[page=7]{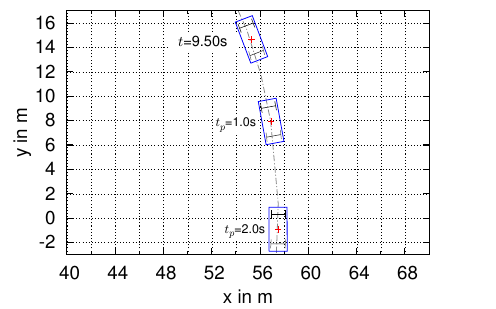}}\label{fig:f5}}
  \hfill
  \subfloat[Deceleration and right turn maneuver.]
           {\resizebox{86mm}{!}{\includegraphics[page=8]{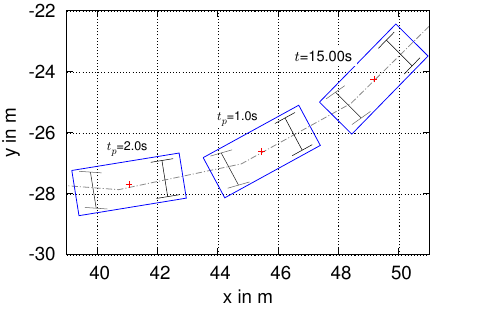}}\label{fig:f6}}
  \caption{Snapshots of the driven course in birds-eye view with two predicted states of the vehicle.
           The cross illustrates the center of gravity.
           For localizing the current state on the reference curve, inspect the coordinates with respect to Figure \ref{fig:lyingEightCurve}.}
  \label{fig:snapshots}
\end{figure*}
The state and input constraints for the \mbox{V-MPC}, are \mbox{${\delta} \in [-20.0^{\circ}, 20.0^{\circ}]$}, \mbox{$u_{\text{1,v}}\in[-5.0, 5.0]^{\circ}/\text{s}$} and \mbox{$u_{\text{2,v}}\in[-2.5, 2.0]\text{m}/\text{s}^2$}.
Defining the desired target state $\boldsymbol{\hat{\xi}}$,  $\hat{v}_{\text{tar}}=10\text{m/s}$ and $\hat{{d_{\perp}}}=0\text{m}$.
\begin{table}[ht]
\caption{Design parameters for both MPC instances.}
\label{MPC parameters}
\begin{center}
\label{tab:params}
\begin{tabular}{ccc}
\toprule
Parameter &\multicolumn{2}{c}{Value}\\
\cmidrule{2-3}
&C-MPC&V-MPC\\
\midrule
$\boldsymbol{Q}  $ & $\text{diag}(1.0, 1.0)$   &  $\text{diag}(0.1, 0.1, 0.2, 0.2, 0.5, 0.5)$\\
$\boldsymbol{R}  $ & $ 0.01 $            &  $\text{diag}(1.0, 0.1)$\\
$\Delta s/\Delta t$ & $0.5\text{m}$      &  $0.05\text{s}$\\
$M       $ & $3$                &  $3$\\
$S_{\text{hor}}/T_{\text{hor}}$ & ${s_{\text{r}}}^{*}_{\text{V-MPC}}(T)$ & $2.0\text{s}$\\
$N_{\text{hor}}   $ & $\text{every 0.5 meters}$             & $20$\\
\bottomrule
\end{tabular}
\end{center}
\end{table}
The C-MPC instance is initialized with the estimated curvature values from equation (\ref{kappaFormula}).
The initial state for the V-MPC is $x_{0,\text{V-MPC}}= [0~0~\pi/4~0~10~0~\pi/4~0]^{\intercal}$.
Figure \ref{fig:stategraphs} depicts the resulting trajectories for the V-MPC states and controls.
The most crucial entity is the lateral displacement, that is ${d_{\perp}}$, which we note remains below one centimeter signed throughout the simulation sequence (see Fig. \ref{fig:f3}).
The velocity profile reflects the consideration of maximum lateral acceleration (see Fig. \ref{fig:f1}).
The captions of Figure \ref{fig:stategraphs} give detailed information.
Figure \ref{fig:snapshots} illustrates two snapshots of the driven route in birds-eye-view.
\subsection{Optimization runtime}
The experiments were conducted on an Intel i7-3520M with 2.9GHz.
Table \ref{tab:rtp} depicts the optimization runtime, that is the MPC optimization time for one time step for each MPC instances.
Clearly, each MPC step of both MPC instances is computed in the submillisecond range, precisely in less than 100$\mu$s, with a mean of 68$\mu\text{s}$ (C-MPC) and 59$\mu\text{s}$ (V-MPC) for all situations along the reference track.
Comparison with a state-of-the-art method, we injected only the vehicle model with equivalent parameters (Table \ref{tab:params}) to the ACADO toolkit \cite{Houska2011a}.
ACADO is a free software package for control and dynamic optimization.
The solution method within the ACADO toolbox relies on a \mbox{SQP-solver} and has a mean computation time for the vehicle model of 209$\mu\text{s}$.
\section{CONCLUSIONS}
We introduced a fast trajectory planning scheme using two MPC instances in a concurrent operation mode, considering nonlinear systems.
This design solves the task both, on a low computational complexity as well as a sufficient precision.
On a theoretical basis, future work will concentrate on stability and convergence statements regarding the concurrent operation scheme approach.
Practically, we will perform experiments on a real test track with an experimental vehicle for automated driving functionality at Ulm University.
\begin{table}[t]
\begin{center}
\caption{Optimization runtime for single MPC step}
\label{tab:rtp}
\begin{tabular}{lccc}
\toprule
& C-MPC & V-MPC & ACADO\\
\cmidrule{2-4}
$t_{\text{CPU}} $ &  68$\mu\text{s}$  &  59$\mu\text{s}$ & 209$\mu\text{s}$\\
\bottomrule
\end{tabular}
\end{center}
\end{table}

\end{document}